\documentclass{article}

\usepackage{cite}
\usepackage{amsmath,amssymb,amsfonts}
\usepackage{algorithmic}
\usepackage[dvipdfm]{graphicx}
\usepackage{bmpsize}
\usepackage{color}
\usepackage{textcomp}
\usepackage{xcolor}
\usepackage{url}
\usepackage{algorithm,algorithmic}
\def\BibTeX{{\rm B\kern-.05em{\sc i\kern-.025em b}\kern-.08em
    T\kern-.1667em\lower.7ex\hbox{E}\kern-.125emX}}
\begin{document}

\title{Performance evaluation of accelerated real and complex multiple-precision sparse matrix--vector multiplication
\thanks{This work was supported by JSPS KAKENHI with Grant Number 23K11127.}
}

\author{Tomonori Kouya\\ORCID: 0000-0003-0178-5519\\hizuoka Institute of Science and Technology}

\maketitle

\begin{abstract}
Sparse matrices have recently played a significant and impactful role in scientific computing, including artificial intelligence-related fields. According to historical studies on sparse matrix--vector multiplication (SpMV), Krylov subspace methods are particularly sensitive to the effects of round-off errors when using floating-point arithmetic. By employing multiple-precision linear computation, convergence can be stabilized by reducing these round-off errors. In this paper, we present the performance of our accelerated SpMV using SIMD instructions, demonstrating its effectiveness through various examples, including Krylov subspace methods.
\end{abstract}


\section{Introduction}

In scientific and technical computing, along with its advanced form, artificial intelligence (AI), the foundation is based on information and communications technology that enables extensive computational processing. Optimized basic linear computations are essential for efficiently utilizing limited computing resources. Accelerating basic linear algebra subprograms (BLAS)\cite{blas} has sternly been a priority with each new hardware generation. These optimizations contribute to higher-level linear algebra libraries such as LAPACK\cite{lapack} and scientific modules in Python, including NumPy, SciPy, TensorFlow, and PyTorch, which are integral to AI development.

As computational scale increases, however, the limitations of traditional IEEE binary32 (24-bit mantissa) and binary64 (53-bit mantissa) BLAS emerge, particularly for ill-conditioned problems where simple application fails to ensure reliability. Currently, an approach that remains relevant involves using multiple-precision floating-point (MPF) arithmetic, which increases mantissa length without altering the algorithm to mitigate round-off errors. Additionally, methods such as the Ozaki scheme allow for higher-speed calculations compared with the MPF arithmetic by utilizing fast binary32 and binary64 BLAS while controlling round-off errors. Therefore, selecting the optimal computation method based on problem conditions and computing environment is essential.

We have previously developed BNCmatmul\cite{bncmatmul}, with multiple-precision well-tuned BLAS and LU decomposition for real and complex dense matrices. A key feature of this development is the use of AVX2 SIMD instructions for faster performance and OpenMP parallelization in shared memory environments, outperforming MPLAPACK/MPBLAS\cite{mplapack}.

This study focuses on the acceleration of multiple-precision sparse matrix--vector multiplication (SpMV). We report the implementation results using AVX2, which has proven performance, along with parallelization using OpenMP. For multiple-precision SpMV, previous studies conducted by Kotakemori et al., summarized in the Lis library \cite{lis}, demonstrate the use of a sparse matrix--vector multiplication that supports the double--double (DD, 106 bits) precision arithmetic, which is a type of multi-component precision. Additionally, the study conducted by Hishinuma et al. \cite{hishinuma2014acs} introduced the block sparse row format for cache optimization, accelerating the performance, with their source code publicly available.

Our SpMV implementation supports not only DD but also triple--double (TD, 159 bits), quad--double (QD, 212 bits), and arbitrary precision using MPFR, covering both real and complex SpMV. This allows us to construct a more comprehensive and high-performance solution than Lis within our BNCmatmul.

This paper begins by explaining the target MPF arithmetic, specifically SIMDization for multi-component types such as DD, TD, and QD precision arithmetic using binary64. Thereafter, we introduce the SpMV algorithm using the CSR format for sparse matrices. Using matrices from the SuiteSparse Matrix Collection \cite{ufsparse}, we evaluate the performance of real and complex SpMV on multi-core CPUs to demonstrate the achieved acceleration. As an application, we also show the acceleration of Krylov subspace methods using real SpMV. Finally, we present our conclusions and future work.

\section{SIMDized multi-component-type multiple-precision arithmetic}

The implemented real SpMV in this study supports multi-component-way formats: DD, TD, QD, and multi-digit-way format: MPFR. It also supports input and output in the coordinate (COO) format using the matrix market (MTX) format, whereas SpMV is implemented as a multiplication of sparse matrices in the compressed sparse row (CSR) format with dense vectors. By incorporating AVX2 into the multi-component precision SpMV, we perform benchmark tests to measure the achievable acceleration, as detailed later.

As shown in the study by Hishinuma et al. \cite{dd_avx_original}, basic linear computations using multi-component MPF arithmetic can be accelerated with SIMD instructions. However, a comprehensive implementation covering TD and QD precision arithmetic appears to be unique to BNCmatmul. Here, we introduce the DD, TD, and QD precision arithmetic used in our implementation, including multiplication with binary64.

In this study, we use the \verb|_m256d| data type, which combines four binary64 values. We utilize the functions \verb|_mm256_[add, sub, mul, div]_pd|, which allow arithmetic operations on this type from C, along with the \verb|_mm256_fmadd_pd| function, which corresponds to fused multiply-add (FMA). Using these functions, we SIMDize the primary functions of error-free transformation techniques---QuickTwoSum, TwoSum, and TwoProd-FMA---and implement them as AVX2QuickTwoSum, AVX2TwoSum, and AVX2TwoProd-FMA functions, respectively.

Below, let $a$, $b$, $c$, and $d$ be of the \verb|_m256d| data type, each representing four binary64 floating-point numbers, such that $a = (a_0, a_1, a_2, a_3)$, $b = (b_0, b_1, b_2, b_3)$, $s = (s_0, s_1, s_2, s_3)$, and $e = (e_0, e_1, e_2, e_3)$.

\paragraph{Addition and multiplication in DD precision}
DD precision arithmetic is constructed as a straightforward combination of these error-free transformation functions; thus, the addition and multiplication operations used for matrix multiplication can be naturally SIMDized and implemented. We have already implemented these as AVX2DDadd (Algorithm \ref{algo:AVX2DDadd}) and AVX2DDmul (Algorithm \ref{algo:AVX2DDmul}) \cite{kouya_iccsa2021}.

\begin{algorithm}[htbp]
    \caption{ $r[2] :=\mathrm{AVX2DDadd}(x[2], y[2])$ }\label{algo:AVX2DDadd}
    \begin{algorithmic}
        \STATE $(s, e) := \mathrm{AVX2TwoSum}(x[0], y[0])$
        \STATE $w := $\verb|_mm256_add_pd|($x[1]$, $y[1]$)
        \STATE $e := $\verb|_mm256_add_pd|($e$, $w$)
        \STATE $(r[0], r[1]) := \mathrm{AVX2QuickTwoSum}(s, e)$
        \STATE \textbf{return} ($r[0]$, $r[1]$)
    \end{algorithmic}
\end{algorithm}

\begin{algorithm}[htbp]
    \caption{ $r[2] :=\mathrm{AVX2DDmul}(x[2], y[2])$ }\label{algo:AVX2DDmul}
    \begin{algorithmic}
        \STATE $(p_1, p_2) := \mathrm{AVX2TwoProd-FMA}(x[0], y[0])$
        \STATE $w_1 := $\verb|_mm256_mul_pd|($x[0], y[1]$)
        \STATE $w_2 := $\verb|_mm256_mul_pd|($x[1], y[0]$)
        \STATE $w_3 := $\verb|_mm256_add_pd|($w_1, w_2$)
        \STATE $p_2 := $\verb|_mm256_add_pd|($p_2, w_3$) 
        \STATE $(r[0], r[1]) :=\mathrm{AVX2QuickTwoSum}(p_1, p_2)$
    \end{algorithmic}
\end{algorithm}

The sparse matrices from the SuiteSparse Matrix Collection used in the benchmarks are all provided in either integer or binary64 precision. Therefore, to maintain SpMV accuracy while conserving memory usage, we implemented binary64 and DD precision floating-point multiplication (DD-D mixed precision multiplication) for SpMV with binary64 sparse matrices as AVX2DDmulD (Algorithm \ref{algo:AVX2DDmulD})

\begin{algorithm}[htbp]
    \caption{ $r[2] :=\mathrm{AVX2DDmulD}(x[2], y)$ }\label{algo:AVX2DDmulD}
    \begin{algorithmic}
        \STATE $(p_1, p_2) := \mathrm{AVX2TwoProd-FMA}(x[0], y)$
        \STATE $p_2 := $\verb|_mm256_add_pd|($x[1], p_2$) 
        \STATE $(r[0], r[1]) :=\mathrm{AVX2QuickTwoSum}(p_1, p_2)$
    \end{algorithmic}
\end{algorithm}

%
\paragraph{Addition and multiplication in QD precision}

For QD arithmetic, we implemented AVX2QDadd and AVX2QDmul using the AVX2-based ThreeSum and ThreeSum2, along with a partially AVX2-optimized Renorm function, based on the sloppy version of addition and multiplication, which involves fewer computations.

We also implemented AVX2QDmulD (Algorithm \ref{algo:AVX2QDmulD}), which performs multiplication between binary64 and QD precision floating-point numbers.

\begin{algorithm}[htbp]
    \caption{ $r[4] :=\mathrm{AVX2QDmulD}(x[4], y)$ }\label{algo:AVX2QDmulD}
    \begin{algorithmic}
        \STATE $(p_0, q_0) := \mathrm{AVX2TwoProd}(x[0], y)$
        \STATE $(p_1, q_1) := \mathrm{AVX2TwoProd}(x[1], y)$
        \STATE $(p_2, q_2) := \mathrm{AVX2TwoProd}(x[2], y)$
        \STATE $p_3 := $\verb|_mm256_add_pd|$(x[3], y)$
        \STATE $s_0 := p_0$
        \STATE $(s_1, s_2) := \mathrm{AVX2TwoProd}(q_0, p_1)$
        \STATE $\mathrm{AVX2ThreeSum}(s_2, q_1, p_2)$   
        \STATE $\mathrm{AVX2ThreeSum2}(q_1, p_2, p_3)$     
        \STATE $s_3 := q_1$
        \STATE $s_4 := $\verb|_mm256_add_pd|($q_2$, $p_2$)
        \STATE $(r[0], r[1], r[2], r[3]) := \mathrm{AVX2Renorm}(s_0, s_1, s_2, s_3, s_4)$
        \STATE \textbf{return} ($r[0], r[1], r[2], r[3]$)
    \end{algorithmic}
\end{algorithm}
%
\paragraph{Addition and multiplication in TD precision}
Optimized triple-precision floating-point arithmetic has been proposed by Fabiano et al. \cite{triple_word_prec2019}; however, acceleration is challenging even with SIMD optimization. Therefore, in this study, we implemented TDaddq and TDmulq by setting $x[3] = y[3] = 0$ for triple precision based on QD precision addition and multiplication and used the SIMD-optimized versions, AVX2TDaddq and AVXTDmulq. We also implemented AVX2TDmulD (Algorithm \ref{algo:AVX2TDmulD}) for multiplication between binary64 and TD precision floating-point numbers.

\begin{algorithm}[htbp]
    \caption{ $r[4] :=\mathrm{AVX2TDmulD}(x[3], y)$ }\label{algo:AVX2TDmulD}
    \begin{algorithmic}
        \STATE $(p_0, q_0) := \mathrm{AVX2TwoProd}(x[0], y)$
        \STATE $(p_1, q_1) := \mathrm{AVX2TwoProd}(x[1], y)$
        \STATE $(p_2, q_2) := \mathrm{AVX2TwoProd}(x[2], y)$
        \STATE $s_0 := p_0$
        \STATE $(s_1, s_2) := \mathrm{AVX2TwoProd}(q_0, p_1)$
        \STATE $\mathrm{AVX2ThreeSum}(s_2, q_1, p_2)$   
        \STATE $s_3 := $\verb|_mm256_add_pd|($q_2$, $p_2$)
        \STATE $(r[0], r[1], r[2]) := \mathrm{AVX2Renorm4}(s_0, s_1, s_2, s_3)$
        \STATE \textbf{return} ($r[0], r[1], r[2], r[3]$)
    \end{algorithmic}
\end{algorithm}

----------------%
\section{CSR Sparse matrix-vector arithmetic and its parallelization}

As previously described, the implemented real SpMV in this study supports four types of MPF arithmetic, including DD, TD, QD, and MPFR. It also supports input and output in the COO format using the MTX format, whereas SpMV is implemented as a multiplication of sparse matrices in the CSR format with dense vectors. For multi-component multiple-precision SpMV, we incorporate AVX2, and, as described later, we perform benchmark tests to measure the degree of acceleration achievable.

For instance, we use the following real square matrix $A\in\mathbb{R}^{5\times 5}$, 
\begin{equation}
A = [a_{ij}]_{i, j = 1, 2, 3, 4, 5} = \left[\begin{array}{ccccc}
	1 & 0 & 2 &  0  &  0 \\
	0 & 3 & 0 &  0  & -4 \\
	0 & 0 & 5 &  0  &  0 \\
	6 & 0 & 0 & -7  &  0 \\
	0 & 0 & 0 &  0  &  8
\end{array}\right] \label{eqn:sparse_sample55}
\end{equation}
In the COO format, which corresponds to the MTX COO format, the elements (data) of the sparse matrix and their positions (row and col) can be specified as follows. The row indices (row) and column indices (col) of the eight elements are stored, respectively.
\begin{verbatim}
data[]   =  [ 1  2  3 -4  5  6 -7  8]
row[]    =  [0 0 1 1 2 3 3 4]
col[]    =  [0 2 1 4 2 0 3 4]
\end{verbatim}

For instance, \( a_{25} = -4 \) is stored as \verb|data[3] = -4|, with \verb|row[3] = 1| (indicating the second row) and \verb|col[3] = 4| (indicating the fifth column). This arrangement, where values are stored in row-wise order (left to right), is known as row-wise format.

For SpMV, the CSR format is used. Elements arranged in row-wise order are stored in \verb|data|, column indices are stored in \verb|indices|, and \verb|indptr| stores the index in \verb|data| of the first element of each row. This reduces the data size compared with the COO format. In this example, for \( a_{25} = -4 \), because \verb|indptr[1] = 2|, the second row starts at the third element, and because \verb|indices[3] = 4|, \verb|data[4] = -4| represents the value of the \( a_{25} \) element.

\begin{verbatim}
data[]   =  [ 1  2  3 -4  5  6 -7  8]
indices[]=  [0 2 1 4 2 0 3 4]
indptr[] =  [0 2 4 5 7 8]
\end{verbatim}

The multiplication of a sparse matrix in CSR format with a dense vector is implemented using AVX2 instructions with 256-bit registers.
\begin{figure}[htb]
    \begin{center}
        \includegraphics[width=.45\textwidth]{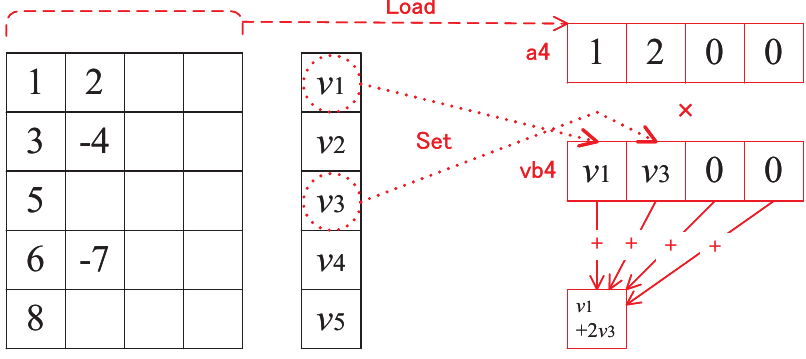}
        \caption{Calculation of SIMDized SpMV}\label{fig:simd_spmv}
    \end{center}
\end{figure}

As shown in \figurename\ \ref{fig:simd_spmv}, elements of the vector corresponding to the non-zero matrix elements in CSR format are extracted, and in the case of AVX2, four elements are arranged together for simultaneous multiplication of the matrix elements with the vector elements, followed by summation. For dense matrix-dense vector multiplication, fast SIMD load/store instructions can be used on contiguous memory regions, but with sparse matrices, vector elements need to be extracted from non-contiguous regions, which can potentially negatively impact the performance. A previous study by Hishinuma et al. \cite{hishinuma2014acs} showed that the structure of a sparse matrix significantly affects cache memory hit rates; thus, actual performance often remains uncertain until execution.

Our multi-component precision matrices and vectors are allocated independently as binary64 matrices and vectors for each component. Therefore, SpMV is executed with two binary64 matrices and vectors for DD precision, three for TD precision, and four for QD precision.

Complex sparse matrices and vectors are represented as real sparse matrices and vectors with the same structure for the real and imaginary parts, respectively. Therefore, for instance, a matrix \( A \in \mathbb{C}^{n \times n} \) and a vector \( \mathbf{v} \in \mathbb{C}^n \) are represented as follows:

\begin{equation}
\begin{split}
    A &= A_\mathrm{re} + \mathrm{i}A_\mathrm{im} \in \mathbb{C}^{n \times n}, \quad A_\mathrm{re}, A_\mathrm{im} \in \mathbb{R}^{n \times n}, \\ 
    \mathbf{v} &= \mathbf{v}_\mathrm{re} + \mathrm{i}\mathbf{v}_\mathrm{im} \in \mathbb{C}^n, \quad \mathbf{v}_\mathrm{re}, \mathbf{v}_\mathrm{im} \in \mathbb{R}^n
\end{split}\label{eqn:complex_expression} 
\end{equation}

and SpMV is executed as

\begin{equation}
    A\mathbf{v} = (A_\mathrm{re}\mathbf{v}_\mathrm{re} - A_\mathrm{im}\mathbf{v}_\mathrm{im}) + \mathrm{i}(A_\mathrm{re}\mathbf{v}_\mathrm{im} + A_\mathrm{im}\mathbf{v}_\mathrm{re})\label{eqn:complex_spmv}
\end{equation}

Thus, the computational cost of complex SpMV is four times that of real SpMV.

\section{Performance evaluation}

This section presents the benchmark test results using the following EPYC environment.
\begin{description}
    \setlength{\parskip}{0cm} 
    \item[H/W] AMD EPYC 9354P 3.7 GHz 32 cores, Ubuntu 20.04.6 LTS
    \item[S/W] Intel Compiler version 2021.10.0, GNU MP 6.2.1, MPFR 4.1.0, MPC 1.2.1, Python 3.8.10
\end{description}
The sparse matrices used were those available in the Sparse Matrix Collection, downloadable via the Python module ssgetpy (\url{https://github.com/drdarshan/ssgetpy}), from which $n \times n$ square sparse matrices $A$ in MTX format were extracted. The $n$-dimensional vector $\mathbf{v}$ is defined as follows:
\begin{equation}
\begin{split}
    \mathrm{Real:} \mathbf{v} &:= \sqrt{2}\ [1\ 2\ ...\ n]^T \\ 
    \mathrm{Complex:} \mathbf{v} &:= \sqrt{2 + 3\mathrm{i}}\ [1\ 2\ ...\ n]^T \\ 
\end{split}\label{eqn:sample_vector}
\end{equation}
Then, the computation time for $A\mathbf{v}$ was measured. Because all matrices $A$ are in binary64 precision, we converted them to the multiple-precision format used for SpMV computations, referred to as P (pure precision arithmetic), whereas using binary64 $A$ for SpMV computations is referred to as M (mixed-precision arithmetic). Arbitrary precision computations using MPFR were executed with a mantissa of 256 bits.

Additionally, for practical purposes, $A^T\mathbf{v}$ is also significant, which is denoted using the abbreviation SpTMV.

\subsection{Real and complex SpMV}

First, to observe the parallelization effects of OpenMP and the performance improvement rate from AVX2, we measured the computation time of SpMV for each sparse matrix using 32 threads and calculated the speedup ratio by dividing it by the time taken with a single-thread execution.


The results for the real SpMV are summarized in \tablename\ \ref{table:real_spmv}. Previous studies have shown that the SpMV performance is highly dependent on matrix size and the structure of the sparse matrix. Therefore, matrix sizes were divided into three groups: \( n = 101 \)-1000, \( n = 1001 \)-5000, and \( n = 5001 \)-10000. For each group, we present the percentage of cases where the speedup ratio exceeded 1, along with the average speedup ratio.

\begin{table}[htb]
\begin{center}
\caption{Speedup ratio of parallelized and SIMDized real SpMV}\label{table:real_spmv}
\begin{tabular}{|c|c|c|c|c|c|}\hline      
\multicolumn{2}{|c|}{32 Threads}& mpfr  &   DD      &	TD     & QD \\ \hline
$n=101$ to $1000$   &P&  9.6\%	        &   51.1\%	& 62.6\%	& 64.4\% \\ 
\# matrices: 219    & &  0.34	        &    29.09	&   5.23	& 5.44 \\ \cline{2-6}
                    &M&   85.8\%	    & 52.5\%	& 56.6\%	& 64.8\% \\
                    & &     2.68	    & 28.06	    & 4.87	    &5.44\\ \hline
$n=1001$ to $5000$  &P&  36.2\%         &	61.2\%	&100.0\%	&100.0\%\\
\# matrices: 309    & &  0.34	        &    29.09	&   5.23	& 5.44 \\ \cline{2-6}
                    &M&   99.7\%	    & 98.4\%	& 99.7\%	& 100.0\% \\
                    & &   1.09	        &6.62	    &7.63	    &11.37                    \\ \hline
$n=5001$ to $10000$ &P&  99.0\%	        &   100.0\%	& 100.0\%	& 100.0\%   \\
\# matrices: 103    & &  2.90	        &   16.89	&   9.44	& 11.31 \\ \cline{2-6}
                    &M&   100.0\%	    &   100.0\%	& 100.0\%	& 100.0\% \\
                    & &   6.32	        &   7.45	&10.84	    &14.32                    \\ \hline
\end{tabular}
\end{center}
\end{table}

For \( n \leq 1000 \), only 52.5\%--85.8\% of cases achieved speedup, whereas for \( n > 5000 \), approximately all cases achieved speedup. For multi-component formats where SIMD operations are effective, the average performance improvement ranged from 7.45 to 11.31 times.

Next, \tablename\ \ref{table:real_sptmv} summarizes the speedup ratios for SpTMV (\( A^T\mathbf{v} \)).

\begin{table}[htb]
\begin{center}
\caption{Speedup ratio of parallelized and SIMDized real SpTMV}\label{table:real_sptmv}
\begin{tabular}{|c|c|c|c|c|c|}\hline      
\multicolumn{2}{|c|}{32 Threads}  &  mpfr             &   DD      &	TD     & QD \\ \hline
$n=101$ to $1000$   &P&  11.9\%	        &   38.8\%	&   16.0\%  &	14.2\% \\
\# matrices: 219    & &  0.74	            &   20.06	&   2.64	&  0.77 \\ \cline{2-6}
                    &M&   18.3\%&	45.7\%  &	15.5\%	& 15.5\% \\
                    & & 0.84	    & 23.70	&   2.51	& 1.51 \\ \hline
$n=1001$ to $5000$  &P&  23.6\%	        &   95.8\%  &	39.5\%  &	35.0\%\\
\# matrices: 309    & &  0.87	            &   2.69    &	1.37	& 1.32\\ \cline{2-6}
                    &M&   34.6\%          &	99.0\%	& 39.2\%	& 35.0\% \\
                    & &   1.40	& 2.84	& 1.35	& 1.33 \\ \hline
$n=5001$ to $10000$ &P&  44.7\%           &	100.0\% &	98.1\%  &	88.3\% \\
\# matrices: 103    & &  1.50	            &   6.05	&   2.80	&   2.62 \\ \cline{2-6}
      &M&   73.8\%	&   100.0\%	    & 99.0\%	& 86.4\% \\
                    & &   2.46	& 6.06	&2.76	&2.63                    \\ \hline
\end{tabular}
\end{center}
\end{table}

For SpTMV, element retrieval from matrix \( A \) is the same as in SpMV; however, because \( A^T\mathbf{v} \) is computed by dividing the matrix across threads, combining the results at the end using a critical section is necessary. This additional computational effort compared with SpMV tends to reduce the effectiveness of parallelization. In practice, even for \( n > 5000 \), perfect speedup has not been achieved for multi-component TD and QD precision formats. Techniques such as identifying the matrix size and structure in advance and performing a transpose before executing SpMV may be effective.


The acceleration rates for complex SpMV and SpTMV are summarized in\ \tablename\ \ref{table:complex_spmv}. In the Sparse Matrix Collection, few complex sparse matrices exist, and when limited to square matrices with $n=101$ to 10000, only the 16 matrices listed in\ \tablename\ \ref{table:complex_spmat} exist.

\begin{table}
\caption{16 complex sparse matrices from Sparse Matrix Collection}\label{table:complex_spmat}
\begin{center}
\begin{tabular}{|c|c|c|}\hline
Name           & $n$   & \# Nonzero \\\hline
\verb|qc324|    & 324	& 26730 \\ 
\verb|young1c|  & 841	& 4089  \\
\verb|young2c|  & 841	& 4089  \\
\verb|young4c|  & 841	& 4089  \\
\verb|dwg961a|  & 961	& 3405  \\
\verb|dwg961b|  & 961	& 10591 \\
\verb|mhd1280a| & 1280	& 47906 \\
\verb|qc2534|   & 2534	& 463360 \\
\verb|conf5_0-4x4-10|	& 3072	& 119808 \\
\verb|conf5_0-4x4-14|	& 3072	& 119808 \\
\verb|conf5_0-4x4-18|	& 3072	& 119808 \\
\verb|conf5_0-4x4-22|	& 3072	& 119808 \\
\verb|conf5_0-4x4-26|	& 3072	& 119808 \\
\verb|conf6_0-4x4-20|	& 3072	& 119808 \\
\verb|conf6_0-4x4-30|	& 3072	& 119808 \\
\verb|mplate| & 5962	& 142190 \\
\verb|aft02| & 8184	& 127762 \\ \hline
\end{tabular}
\end{center}
\end{table}

Complex SpMV and SpTMV require four times the computational cost of real SpMV and SpTMV, as shown in (\ref{eqn:complex_spmv}). Therefore, the acceleration rate tends to be higher even for matrices of the same size.
\begin{table}[htb]
\begin{center}
    \caption{Speedup ratio of parallelized and SIMDized complex SpMV and SpTMV}\label{table:complex_spmv}
    \begin{tabular}{|c|c|c|c|c|c|}\hline      
    \multicolumn{2}{|c|}{32 Threads}          &  mpfr             &   DD      &	TD     & QD \\ \hline
    SpMV                &P &   70.6\%          &	100.0\%	& 100.0\%	& 100.0\% \\
    \# matrices: 16     & &  5.45	&7.88	&8.30	&10.99    \\ \cline{2-6}
                        &M&  100.0\%	&94.1\%	&100.0\%	& 100.0\% \\
                        & &  13.35	&3.64	&9.74	&14.15                        \\ \hline
    SpTMV               &P&  76.5\%           &	76.5\% &	76.5\%  &	76.5\% \\
    \# matrices: 16     & & 3.49&	4.37&	2.79&	2.79    \\ \cline{2-6}
                        &M&   82.4\%          &	88.2\%	& 76.5\%	& 76.5\% \\
                        & & 4.71&	4.51&	2.66&	2.75                        \\ \hline
                                                                 
\end{tabular}
\end{center}
\end{table}
\subsection{Krylov subspace methods}

Here, as an application of SpMV, we present the results of solving the system of linear equations 
\begin{equation}
    A\mathbf{v} = \mathbf{b}, A\in\mathbb{R}^{n\times n}, \mathbf{b}\in\mathbb{R}^n \label{eqn:linear_eq}
\end{equation}
using Krylov subspace methods. The vector $\mathbf{v}$ is provided in (\ref{eqn:sample_vector}).

\begin{figure}[htb]
    \begin{center}
        \includegraphics[width=.224\textwidth]{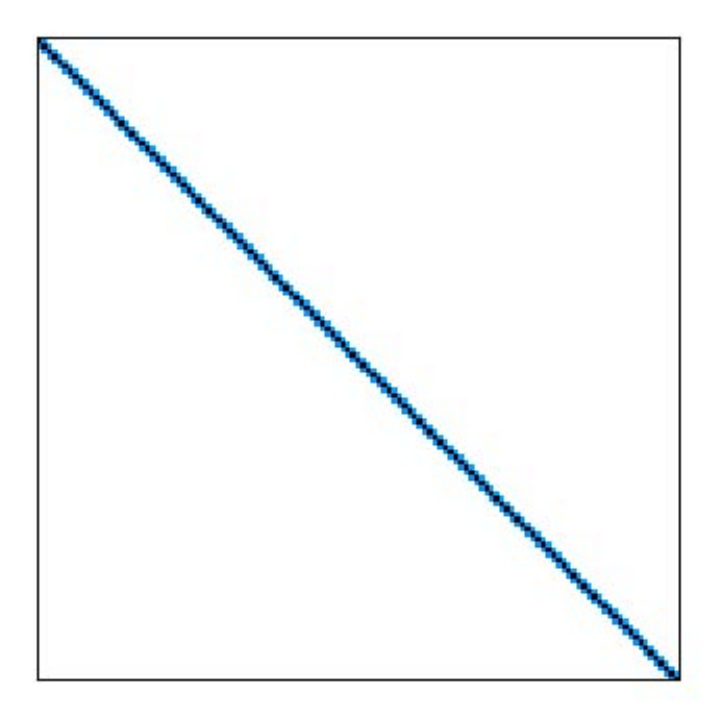}
        \includegraphics[width=.224\textwidth]{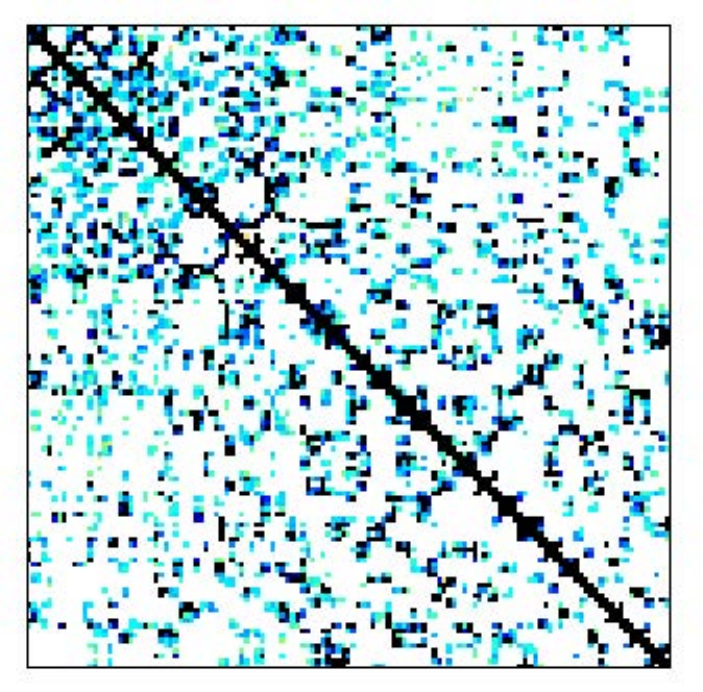}
    \end{center}
    \caption{Structure of tub1000 (left) and nd3k (right)\cite{ufsparse}}\label{fig:tub1000_structure}
\end{figure}

The sparse matrices used were tub1000 (a nonsymmetric matrix) and nd3k (a symmetric positive-definite matrix). Their mathematical properties are listed in \tablename\ \ref{table:tub1000_property}, and the structure of the matrices is shown in \figurename\ref{fig:tub1000_structure}.
\begin{table}[htb]
    \begin{center}
    \caption{Property of tub1000 (left) and nd3k (right)\cite{ufsparse}}\label{table:tub1000_property}
	\begin{tabular}{|c|c|c|c|c|c|}\hline
		        & $n$ & Nonzeros & Symmetric & Cond($A$) \\ \hline
		tub1000 & 1000 &    3996 & No & $1.3\times 10^6$ \\ \hline
		nd3k    & 9000 & 3279690 & Yes & $1.6\times 10^7$ \\ \hline
    \end{tabular}
    \end{center}
\end{table}

For comparison with binary64, relatively well-conditioned matrices were used. Performance evaluations of SpMV and SpTMV for both matrices revealed that the acceleration for tub1000 was not successful, making it meaningful to compare it with nd3k as a benchmark for evaluating the performance of Krylov subspace methods that extensively use SpMV and SpTMV.

\paragraph{Product-type Krylov subspace methods for tub1000}

\begin{figure}[htb]
    \begin{center}
        \includegraphics[width=.45\textwidth]{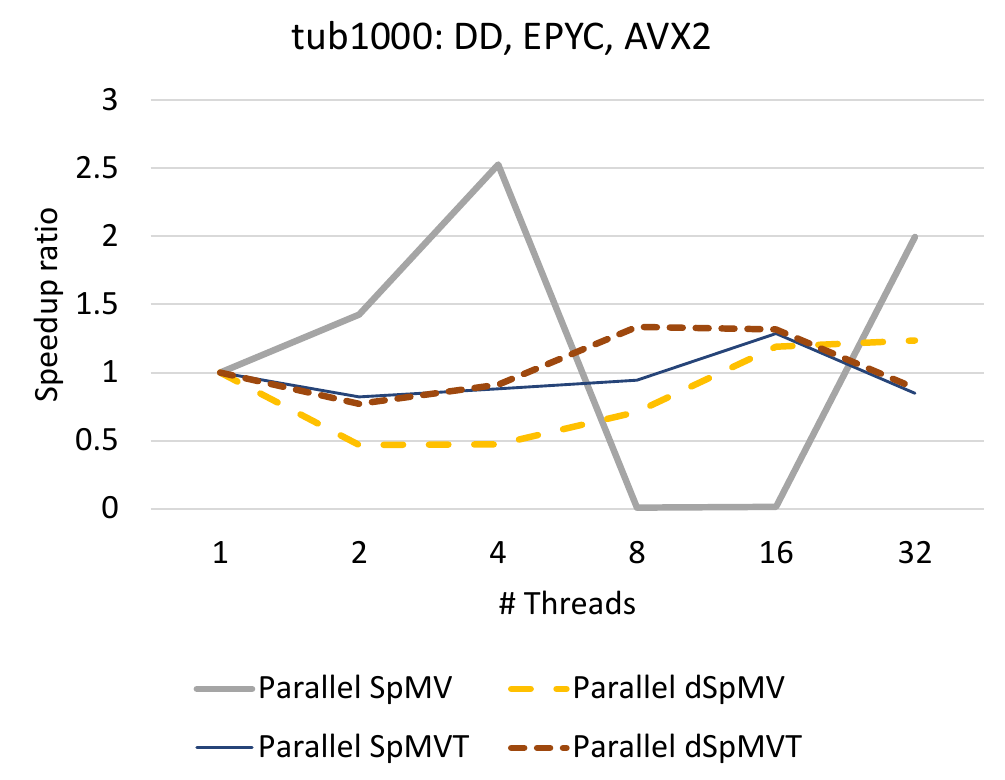}
        \includegraphics[width=.45\textwidth]{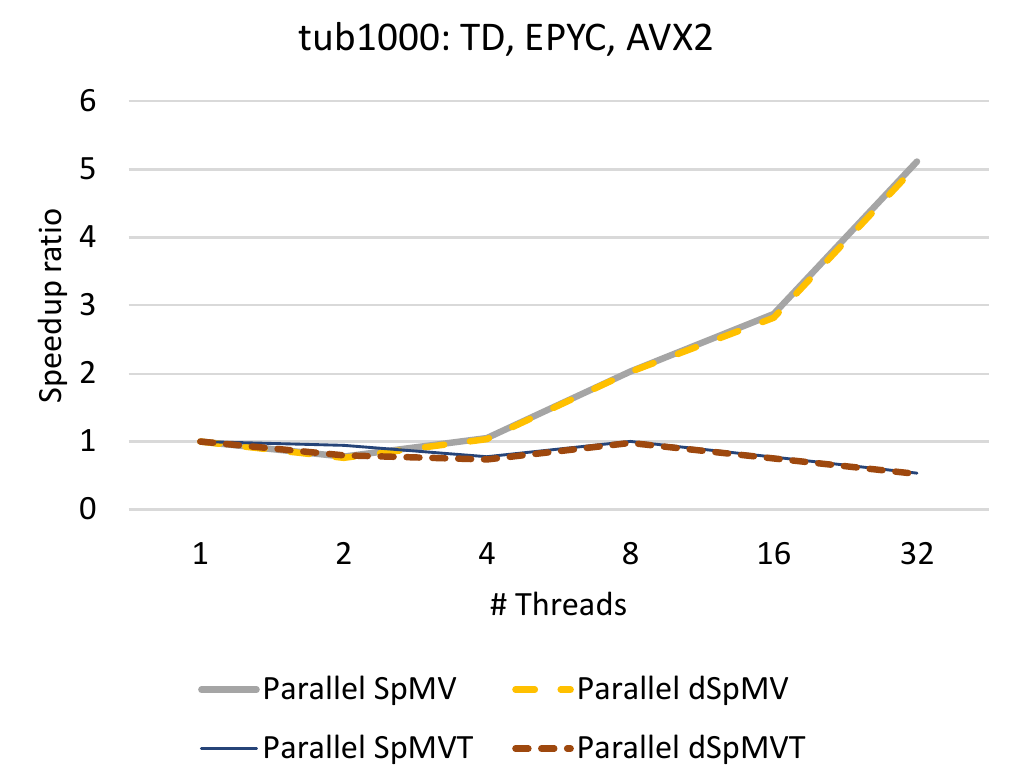}
        \includegraphics[width=.45\textwidth]{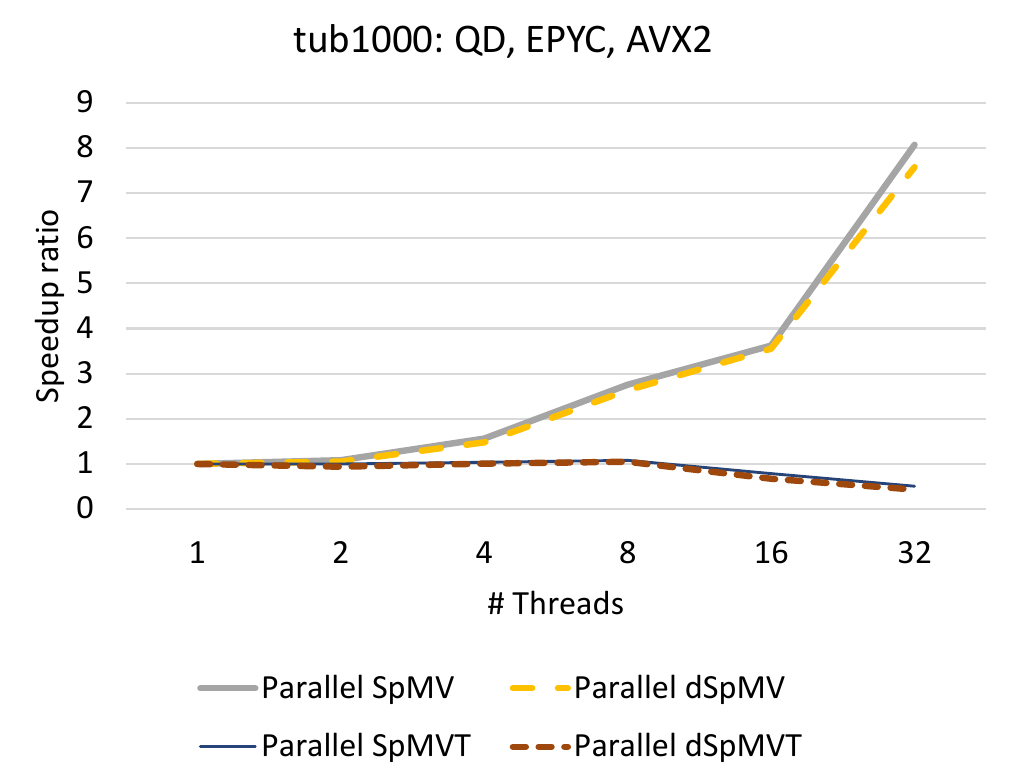}
        \includegraphics[width=.45\textwidth]{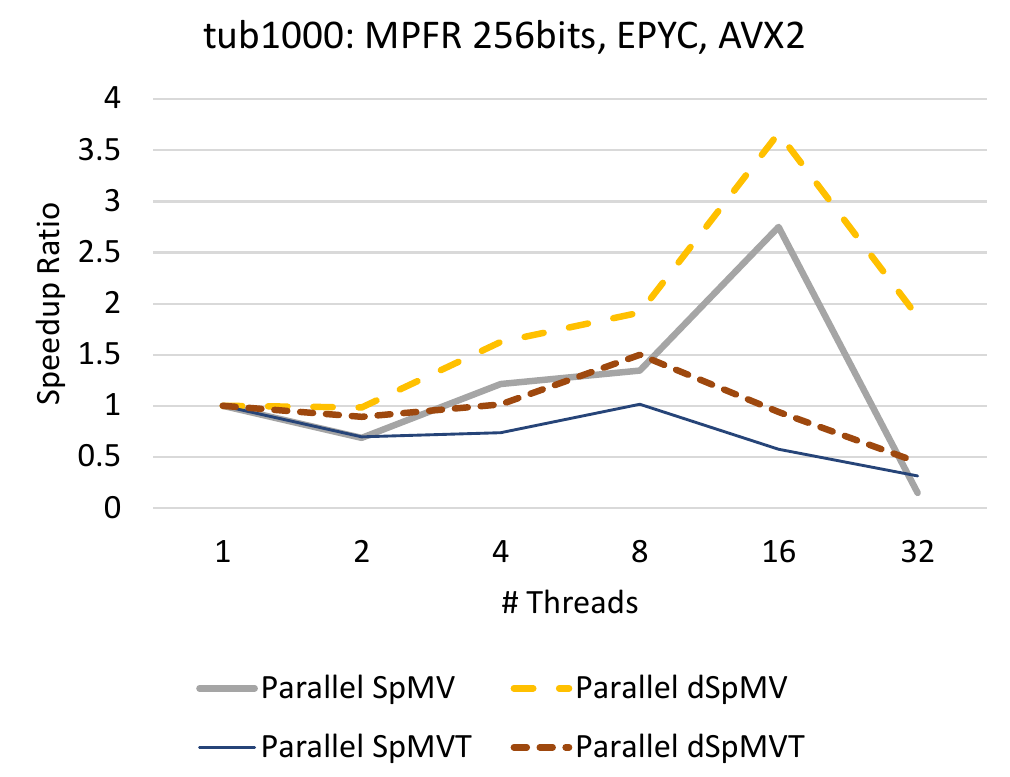}
    \end{center}
    \caption{Speedup ratio of parallelized SpMV, mixed-precision SpMV (dSpMV) for tub1000}\label{fig:spmv_tub1000}
\end{figure}

As product-type Krylov subspace methods applicable to nonsymmetric matrices, we used the BiCG, CGS, BiCGSTAB, and GPBiCG methods to solve the system of linear equations (\ref{eqn:linear_eq}). \figurename\ \ref{fig:spmv_tub1000} shows a plot of the speedup ratios for SpMV and SpTMV when using tub1000 as matrix $A$.

The structure of the matrix closely resembles that of a band matrix, with a narrow bandwidth, which makes it evident that performance improvement, particularly for SpTMV, has not been achieved at all. Furthermore, the use of AVX2 appears to have hindered performance improvements in DD precision SpMV.

\begin{figure}[htb]
    \begin{center}
        \includegraphics[width=.475\textwidth]{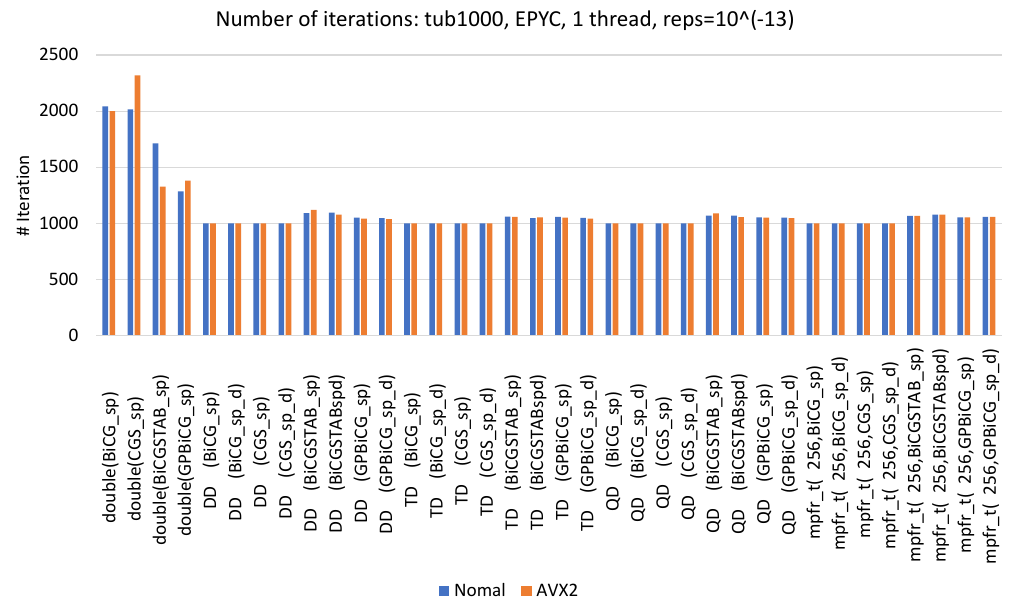}
        \includegraphics[width=.475\textwidth]{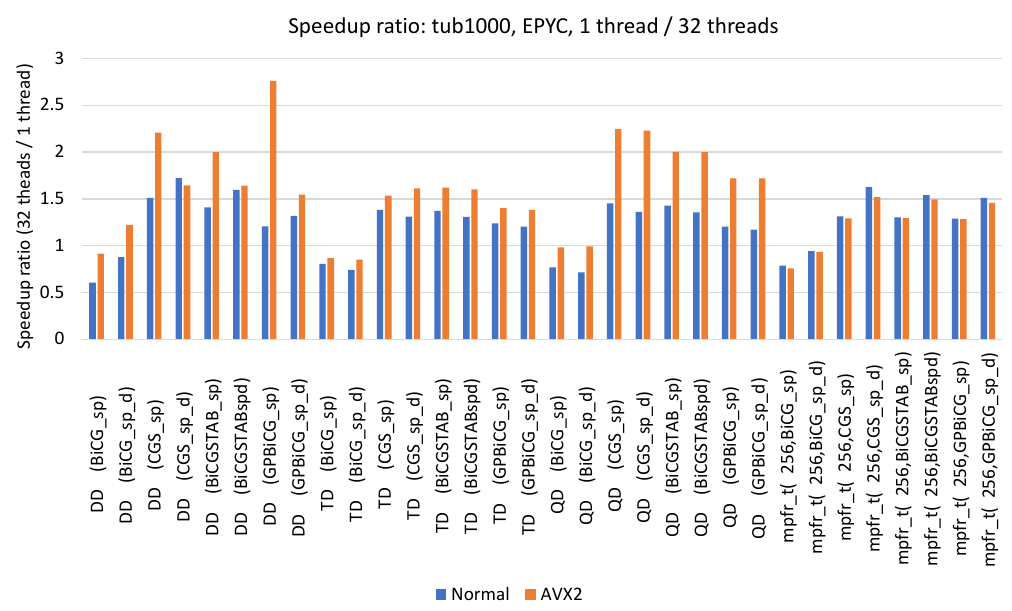}
    \end{center}
    \caption{Number of iterations of Krylov subspace methods (upper) and speedup ratio (lower) for tub1000}\label{fig:num_iter_krylov_speedup_tub1000}
\end{figure}

\figurename\ref{fig:num_iter_krylov_speedup_tub1000} shows the results of solving the system of linear equations with coefficient matrix $A = \text{tub1000}$ using the four methods. The initial guess was set to $\mathbf{x}_0 := 0$, and iterations were terminated when the residual $\mathbf{r_k} := \mathbf{b} - A\mathbf{x}_k$ satisfied $\|\mathbf{r}_k\|_2 \leq 10^{-13}\|\mathbf{b}\|_2 + 10^{-99}$.

As the precision increases, the number of iterations stabilizes at a certain value; however, the speedup ratio varies depending on the computational load of the algorithm. The performance of the BiCG method deteriorates with DD and MPFR precision because of the inclusion of SpTMV. For the other Krylov subspace methods, both TD and QD precision achieve faster speeds, though the improvement ratio is not particularly high.

\paragraph{Conjugate-gradient method for nd3k}

In contrast to tub1000, the SpMV of nd3k, which has a uniformly distributed pattern of non-zero matrix elements, exhibits a high rate of performance improvement through parallelization and SIMD optimization. Because nd3k is a symmetric matrix, the performance difference between SpMV and SpTMV is due to differences in the parallelization approach.

\begin{figure}[htb]
    \begin{center}
        \includegraphics[width=.45\textwidth]{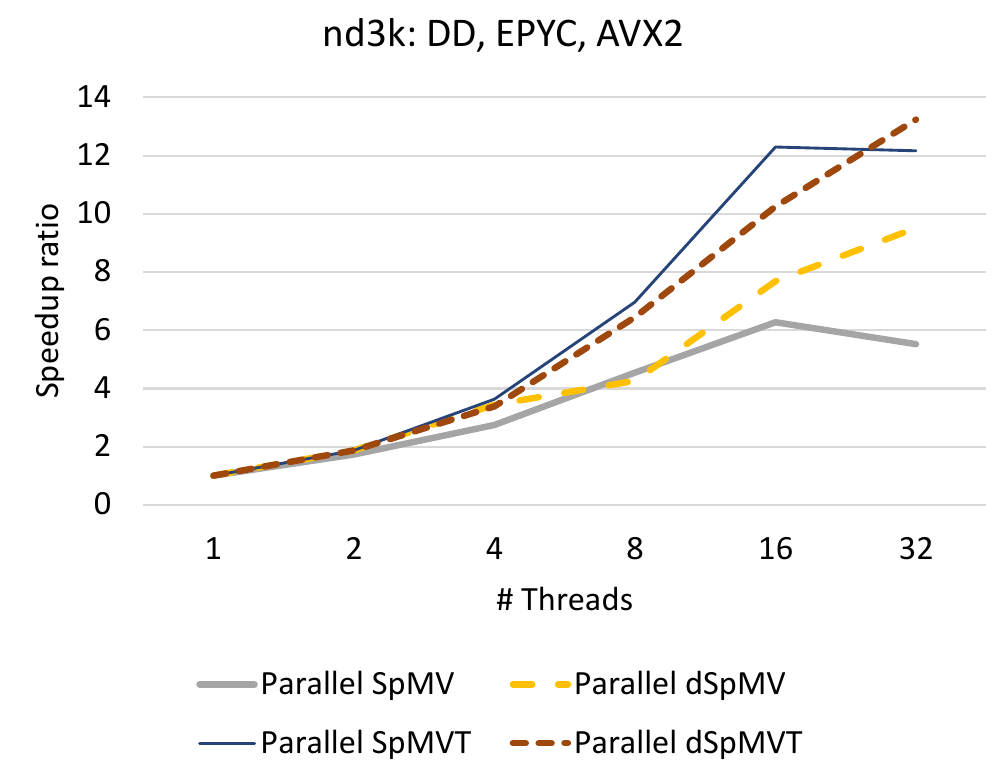}
        \includegraphics[width=.45\textwidth]{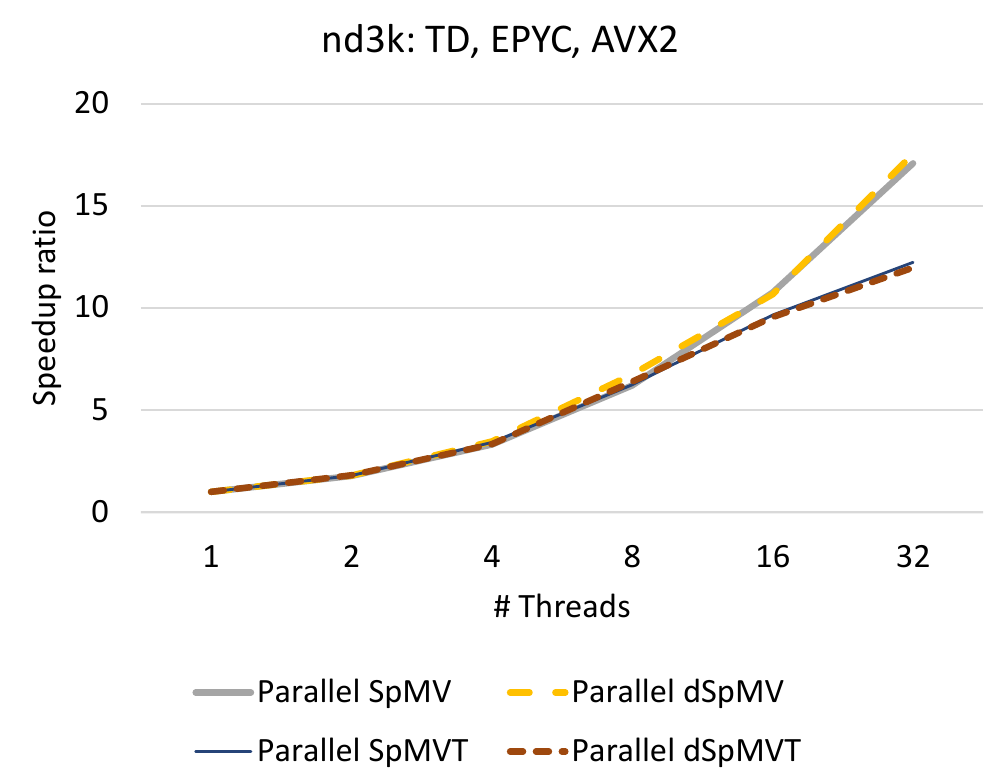}
        \includegraphics[width=.45\textwidth]{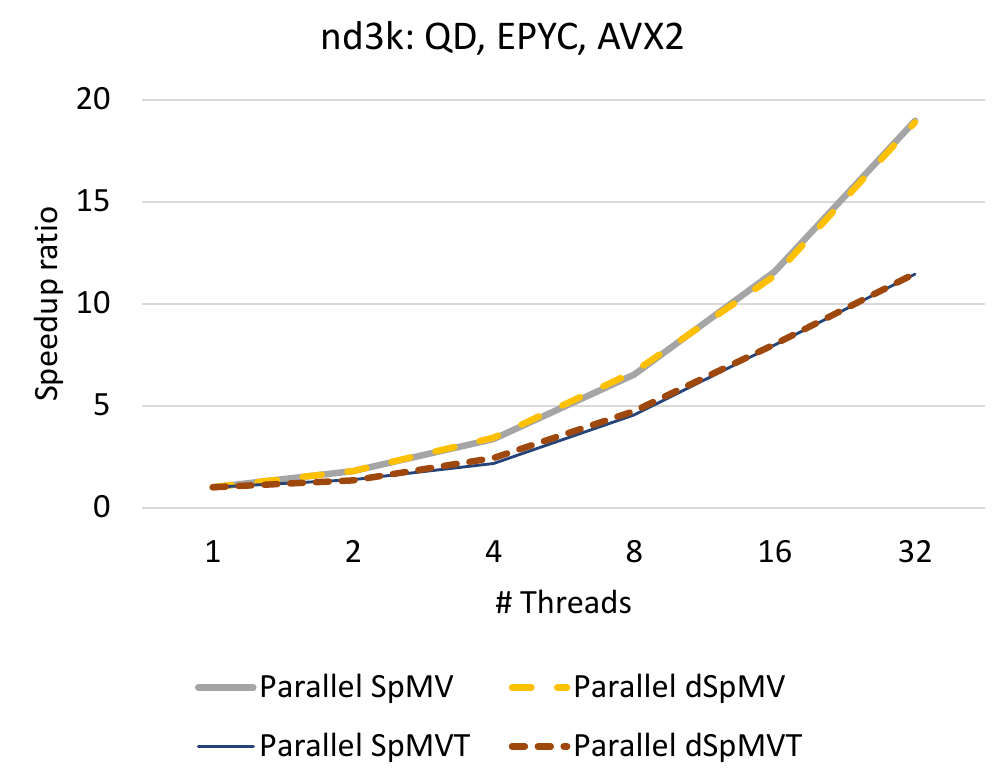}
        \includegraphics[width=.45\textwidth]{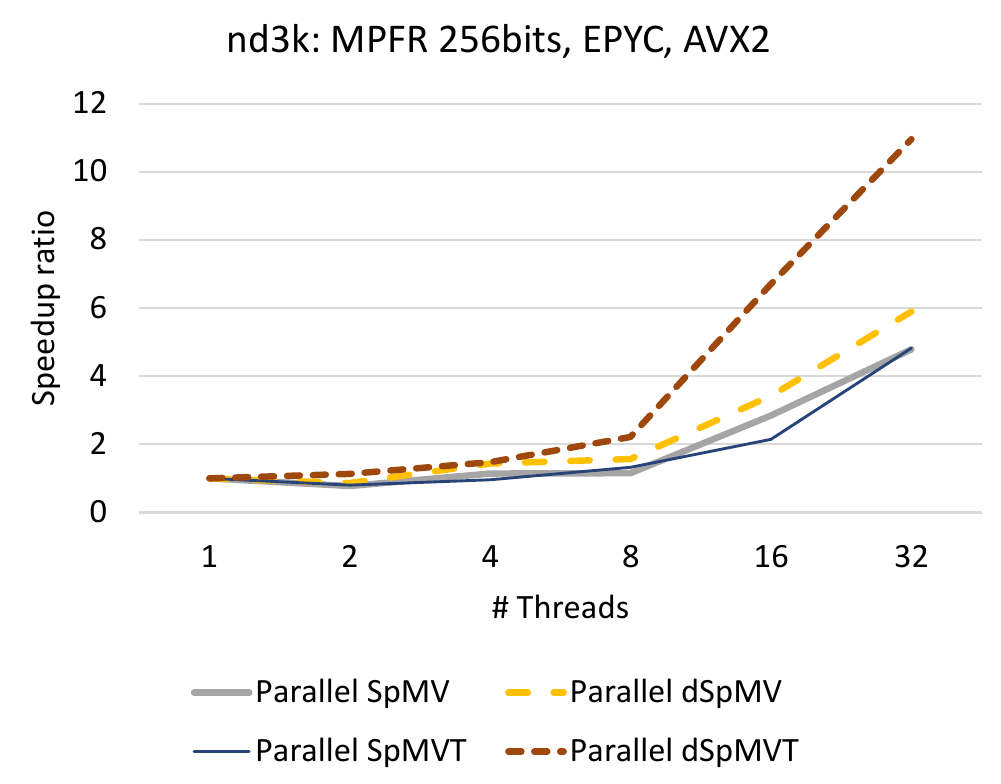}
    \end{center}
    \caption{Speedup ratio of parallelized SpMV, dSpMV for nd3k}\label{fig:spmv_nd3k}
\end{figure}

Performance improvements are achieved at all precision levels. Additionally, SpTMV generally exhibits better performance.

\begin{figure}[htb]
    \begin{center}
        \includegraphics[width=.35\textwidth]{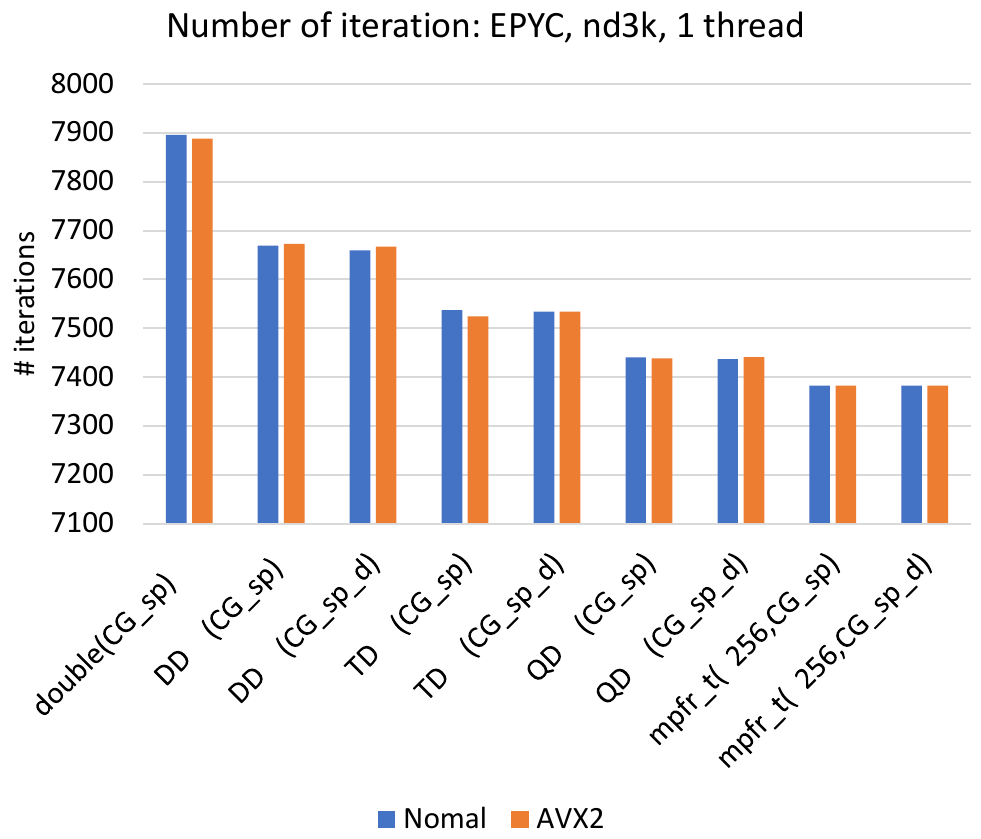}
        \includegraphics[width=.35\textwidth]{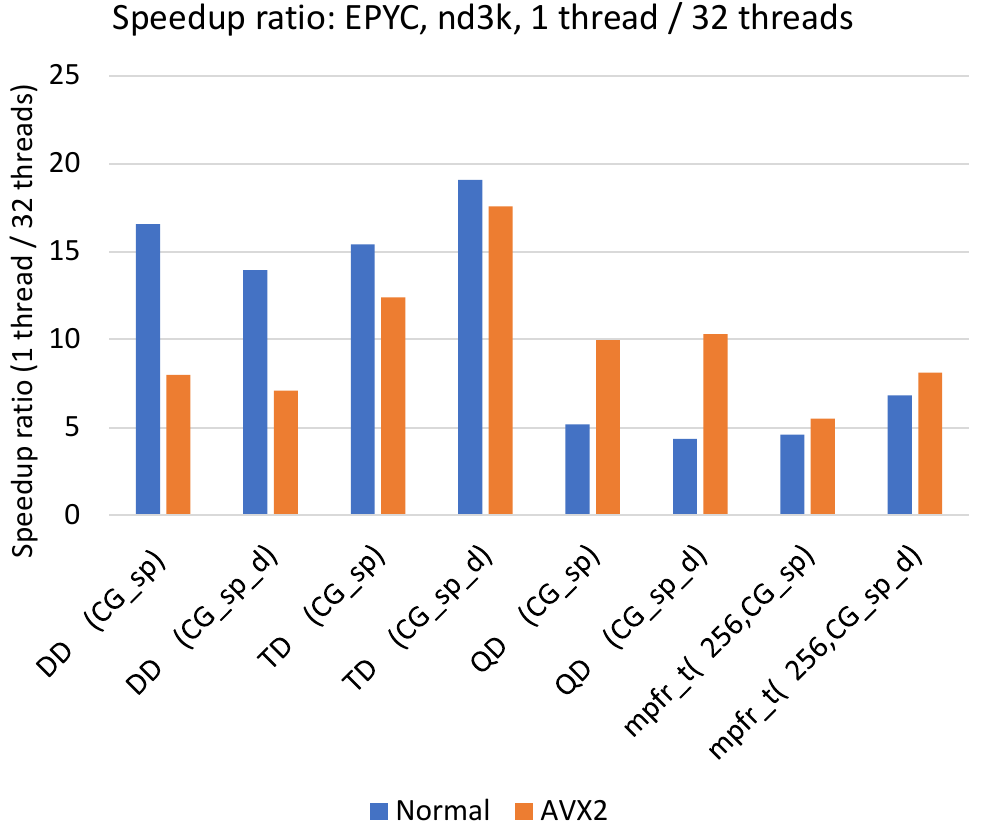}
    \end{center}
    \caption{Number of iterations of CG method and speedup ratio for nd3k}\label{fig:num_iter_cg_speedup_nd3k}
\end{figure}

Because the conjugate-gradient (CG) method can be applied to symmetric positive-definite matrices, we present the performance evaluation results in \figurename\ref{fig:num_iter_cg_speedup_nd3k}. The stopping criterion is the same as that for tub1000.

The number of iterations decreases inversely with precision, and the overall speedup ratio is high.

\section{Conclusion and future work}

The implementation of multiple-precision SpMV confirmed that AVX2 acceleration was achieved in numerous cases for the multi-component DD, TD, and QD precision formats. This acceleration was particularly effective for larger matrix sizes, and it was also confirmed that OpenMP contributed to performance improvements in both real and complex SpMV.

Future studies include exploring practical accuracy enhancement methods for the Lanczos method using multiple-precision calculations and implementing and evaluating Krylov subspace methods that incorporate mixed-precision preconditioning. Additionally, we plan to expand and optimize the functionality of BNCmatmul and extend support for multi-component precision SpMV.

\section*{Acknowledgment}
This work was supported by JSPS KAKENHI with Grant Number 23K11127.



\begin{thebibliography}{10}
    
    \bibitem{blas}
    (2024, Dec.) BLAS. [Online] Available: \url{https://www.netlib.org/blas/}
    
    \bibitem{ufsparse}
    (2024, Dec.) SuiteSparse Matrix Collection. [Online] Available: \url{https://sparse.tamu.edu/}
    
    \bibitem{triple_word_prec2019}
    N.~{Fabiano}, J.-.~M. {Muller}, and J.~{Picot}, ``Algorithms for triple-word
      arithmetic,'' \emph{IEEE Trans. on Computers}, vol.~68, pp. 1573--1583, 2019.
    
    \bibitem{dd_avx_original}
    T.~Hishinuma, A.~Fujii, T.~Tanaka, and H.~Hasegawa, ``Avx acceleration of dd
      arithmetic between a sparse matrix and vector,'' in \emph{Parallel Processing
      and Applied Mathematics}, R.~Wyrzykowski, J.~Dongarra, K.~Karczewski, and
      J.~Wa{\'{s}}niewski, Eds.\hskip 1em plus 0.5em minus 0.4em\relax Berlin,
      Heidelberg: Springer Berlin Heidelberg, 2014, pp. 622--631.
    
    \bibitem{lis}
    (2024, Dec.) T.~Kotakemori, S.~Fujii, H.~Hasegawa, and A.~Nishida, ``Lis: Library of
      iterative solvers for linear systems,'' [Online] Available: \url{https://www.ssisc.org/lis/}.
    
    \bibitem{bncmatmul}
    (2024, Dec.) T.~Kouya, {BNC}matmul. [Online] Available: \url{https://github.com/tkouya/bncmatmul}
    
    \bibitem{kouya_iccsa2021}
    ------, ``Acceleration of multiple precision matrix multiplication based on
      multi-component floating-point arithmetic using avx2,'' in
      \emph{Computational Science and Its Applications -- ICCSA 2021}, O.~Gervasi,
      B.~Murgante, S.~Misra, C.~Garau, I.~Ble{\v{c}}i{\'{c}}, D.~Taniar, B.~O.
      Apduhan, A.~M. A.~C. Rocha, E.~Tarantino, and C.~M. Torre, Eds.\hskip 1em
      plus 0.5em minus 0.4em\relax Cham: Springer International Publishing, 2021,
      pp. 202--217.
    
    \bibitem{lapack}
    (2024, Dec.) LAPACK. [Online] Available: \url{https://www.netlib.org/lapack/}
    
    \bibitem{mplapack}
    (2024, Dec.) M.Nakata, ``Multiple precision arithmetic {LAPACK} and {BLAS}. [Online] Available: \url{https://github.com/nakatamaho/mplapack}
    
    \bibitem{hishinuma2014acs}
    T.~Hishinuma, A.~Fujii, T.~Tanaka, and H.~Hasegawa, ``{AVX2 Acceleration of Double Precision Sparse Matrix in BCRS Format and DD Vector Product},''
      \emph{{IPSJ ACS)}}, vol.~7, no.~4,
      pp. 25--33, 2014. (in Japanese)
    
    \end{thebibliography}

\end{document}